\documentclass[english,10pt]{article}
\usepackage[T1]{fontenc}
\usepackage[latin9]{inputenc}
\usepackage{textcomp}
\usepackage{amsthm}
\usepackage{amsmath}
\usepackage{amssymb}
\theoremstyle{plain}

\usepackage{mathrsfs}
\usepackage{hyperref}
\usepackage{xy}
\hypersetup{colorlinks=true,citecolor=blue}
\usepackage{babel}

\begin{document}

\title{Minimal modularity lifting for GL$_{\text{2}}$ over an arbitrary
number field}

\author{David Hansen%
\thanks{Department of Mathematics, Boston College, Chestnut Hill MA; dvd.hnsn@gmail.com%
}}
\maketitle
\begin{abstract}
We prove a modularity lifting theorem for minimally ramified deformations
of two-dimensional odd Galois representations, over an arbitrary number
field. The main ingredient is a generalization of the Taylor-Wiles
method in which we patch complexes rather than modules.
\end{abstract}

\section{Introduction}

Fix a number field $F/\mathbf{Q}$. The Taylor-Wiles method \cite{TW}
is a technique for proving that a surjection $R_{\overline{\rho}}\twoheadrightarrow\mathbf{T}_{\mathfrak{m}}$
from a Galois deformation ring to a Hecke algebra at minimal level
is an isomorphism. $ $Essentially all incarnations of the Taylor-Wiles
method have been limited to situations where $F$ is totally real
or CM, $R_{\overline{\rho}}$ parametrizes deformations satisfying
strong regularity and self-duality assumptions, and $\mathbf{T}_{\mathfrak{m}}$
arises from the the middle-dimensional cohomology of a Shimura variety.
In a recent and very striking breakthrough, Calegari and Geraghty
 \cite{CG} found a novel generalization of the Taylor-Wiles method
which eliminates some of these restrictions. More precisely, their
method applies when $\mathbf{T}_{\mathfrak{m}}$ acts on the cohomology
of a locally symmetric space $X$ such that $H^{i}(X)_{\mathfrak{m}}$
is nonvanishing in only two consecutive degrees. In this paper we
develop a further generalization of the Taylor-Wiles method; in principle,
our method requires no restriction whatsoever on the range of degrees
for which $H^{\ast}(X)_{\mathfrak{m}}$ is nonzero.

As a sample application, we prove the following theorem, restricting
ourselves to the simplest possible situation in which our technique
yields a new result. Let $F$ be an \emph{arbitrary }number field;
set $d=[F:\mathbf{Q}]$, and let $r$ be the number of nonreal infinite
places of $F$. Fix a finite field $k$ of characteristic $p\geq3$ with $p$ unramified in $F$,
and set $\mathcal{O}=W(k)$. Fix an absolutely irreducible Galois
representation $\overline{\rho}:G_{F}\to\mathrm{GL}_{2}(k)$ unramified
at all but finitely many primes, with no {}``vexing'' primes of
ramification. Suppose $\overline{\rho}$ has the following properties:
\begin{itemize}
\item $\overline{\rho}|D_{v}$ is ordinary or finite flat for all $v|p$,
\item $\det\overline{\rho}(c_{\sigma})=-1$ for all real infinite places
$\sigma$ and complex conjugations $c_{\sigma}$,
\item $\overline{\rho}|G_{F(\zeta_{p})}$ is absolutely irreducible.
\end{itemize}
Let $\mathcal{D}$ denote the functor which assigns to an Artinian
$\mathcal{O}$-algebra $A$ the set of equivalence classes of deformations
$\tilde{\rho}:G_{F}\to\mathrm{GL}_{2}(A)$ of $\overline{\rho}$ which
are minimally ramified at all places $v\nmid p$ and ordinary or finite
flat at all places $v|p$. This functor is represented by a complete
local Noetherian $\mathcal{O}$-algebra $R_{\overline{\rho}}^{\mathrm{min}}$
together with a natural universal lifting $\rho^{\mathrm{min}}:G_{F}\to\mathrm{GL}_{2}(R_{\overline{\rho}}^{\mathrm{min}})$.
Let $\mathbf{T}$ be the Hecke algebra defined in §3; this is defined
as a subalgebra of the ring of endomorphisms of $H^{\ast}(Y,\mathcal{O})$
for $Y$ a certain locally symmetric quotient of $\mathrm{GL}_{2}(F_{\infty})$.
We suppose there is a maximal ideal $\mathfrak{m}\subset\mathbf{T}$
with residue field $k$ together with a surjection $\phi_{\mathfrak{m}}:R_{\overline{\rho}}^{\mathrm{min}}\twoheadrightarrow\mathbf{T}_{\mathfrak{m}}$
such that $\rho_{\mathfrak{m}}=\phi_{\mathfrak{m}}\circ\rho^{\mathrm{min}}:G_{F}\to\mathrm{GL}_{2}(\mathbf{T}_{\mathfrak{m}})$
has characteristic polynomial $X^{2}-T_{v}X+\mathbf{N}v\left\langle v\right\rangle \in\mathbf{T}_{\mathfrak{m}}[X]$
on almost all Frobenius elements $\mathrm{Frob}_{v}$. In order to
apply our generalization of the Taylor-Wiles method we also need this
in non-minimal situations, and we need to know something about the
vanishing of cohomology after localizing at a {}``non-Eisenstein''
prime; for a precise statement, see Conjecture 3.1.

\textbf{Theorem 1.1. }\emph{Suppose Conjecture 3.1 is true. Then $ $$\phi_{\mathfrak{m}}:R_{\overline{\rho}}^{\mathrm{min}}\to\mathbf{T}_{\mathfrak{m}}$
is an isomorphism and $H^{d}(Y,\mathcal{O})_{\mathfrak{m}}$ is free
over $\mathbf{T}_{\mathfrak{m}}$.}

When $r=1$ this theorem follows from the method of  \cite{CG}. Note
that $\mathbf{T}_{\mathfrak{m}}$ often contains $\mathcal{O}$-torsion
elements when $r\geq1$, whereas the classical Taylor-Wiles method
(as streamlined by Diamond  \cite{Diamond} and Fujiwara  \cite{Fuj})
requires an \emph{a priori }assumption that $\mathbf{T}_{\mathfrak{m}}$
be $\mathcal{O}$-flat. 

Let us briefly explain the proof of Theorem 1.1. Set $q=\mathrm{dim}_{k}H_{\emptyset}^{1}(F,\mathrm{ad}^{0}\overline{\rho}(1))$,
and write $R_{\infty}=\mathcal{O}[[x_{1},\dots,x_{q-r}]]$ and $S_{\infty}=\mathcal{O}[[T_1,\dots,T_q]]$.  Let us abbreviate $H=H^{d}(Y,\mathcal{O})_{\mathfrak{m}}$ and $R=R_{\overline{\rho}}^{\mathrm{min}}$; we regard $H$ as an $R$-module via $\phi_{\mathfrak{m}}$. By a patching technique (Theorem 2.2.1), we construct an algebra homomorphism
$i_{\infty}:S_{\infty}\to R_{\infty}$ and a finite $R_{\infty}$-module
$H_{\infty}$, together with a surjection $ $$\phi_{\infty}:R_{\infty}\twoheadrightarrow R$
and an ideal $\mathfrak{a}\subset S_{\infty}$ with $(\phi_{\infty}\circ i_{\infty})(\mathfrak{a})=0$
$ $such that $H\simeq H_{\infty}/\mathfrak{a}H_{\infty}$ as $R_{\infty}$-modules,
where $R_{\infty}$ acts on $H$ through $\phi_{\infty}$. Suppose
we could show the $S_{\infty}$-depth of $H_{\infty}$ was at least
$1+q-r$. Then via $i_{\infty}$ the $R_{\infty}$-depth would be
at least $1+q-r=\mathrm{dim}R_{\infty}$, so $H_{\infty}$ would be
free over $R_{\infty}$ by the Auslander-Buchsbaum formula. We would
then easily conclude that $H$ is free over $R_{\infty}/i_{\infty}(\mathfrak{a})$,
whence the surjection $R_{\infty}/i_{\infty}(\mathfrak{a})\twoheadrightarrow R$
would be an isomorphism and $H$ would be free over $R$. 

In order to carry this out, we appeal crucially to the construction
of $H_{\infty}$: it is the top degree cohomology of a complex $F_{\infty}^{\bullet}$
of free finite rank $S_{\infty}$-modules concentrated in a range
of degrees of length $\leq r$. By a general theorem in commutative
algebra (Theorem 2.1.1), this forces every irreducible component of
the $S_{\infty}$-support of $H^{\ast}(F_{\infty}^{\bullet})$ to
have dimension $\geq1+q-r$. However, the patching construction yields
an $R_{\infty}$-module structure on $H_{\infty}^{\ast}$ which implies
the opposite inequality, from whence we deduce (by Theorem 2.1.1 again)
that $H^{i}(F_{\infty}^{\bullet})$ vanishes for all degrees $i$
except the top degree. As such, $F_{\infty}^{\bullet}$ yields a free
resolution of $H_{\infty}$ of length $r$, so $\mathrm{projdim}_{S_{\infty}}(H_{\infty})=r$.
But then $\mathrm{depth}_{S_{\infty}}(H_{\infty})=1+q-r$ by another
application of Auslander-Buchsbaum.

The numerical coincidence driving this argument persists far beyond
$\mathrm{GL}_{2}$. Roughly speaking, when considering a Galois representation
$\overline{\rho}:\mathrm{Gal}(\overline{F}/F)\to\widehat{G}(k)$ for
$G$ some ($F$-split) reductive algebraic group, we require the equality\[
[F:\mathbf{Q}](\mathrm{dim}G-\mathrm{dim}B)+l(G)=\sum_{v|\infty}H^{0}(F_{v},\mathrm{ad}^{0}\overline{\rho})\]
where $l(G)$ denotes the length of the range of degrees for which
deformations of $\overline{\rho}$ contribute to the Betti cohomology
of locally symmetric quotients of $G(F_{\infty})$; the reader may
wish to compare this with the numerical condition given in  \cite{CHT}.
At the very least, our method generalizes to the case when $\overline{\rho}:G_{F}\to\mathrm{GL}_{n}(k)$
is odd (i.e. $|\mathrm{tr}\overline{\rho}(c_{\sigma})|\leq1$ for
all real places $\sigma$ and complex conjugations $c_{\sigma}$)
and absolutely irreducible with big image, and $\mathcal{D}$ parametrizes
minimally ramified regular crystalline deformations in the Fontaine-Laffaille
range. Note the absence of any restrictions on $F$ or any self-duality
hypothesis on $\overline{\rho}$. This, again, is contingent on assuming
the existence of various surjections $R_{\overline{\rho}}\twoheadrightarrow\mathbf{T}_{\mathfrak{m}}$.

\subsection*{Acknowledgments}

The great debt of inspiration this work owes to the beautiful ideas
of Calegari and Geraghty  \cite{CG} will be evident to the reader;
the idea of patching complexes grew naturally out of their success
at patching presentations. In addition, I'm grateful to Calegari and Geraghty for pointing out a mistake in the initial public version of this paper. I'm also grateful to Avner Ash, Hailong Dao, Michael Harris, and Jack Thorne for some helpful remarks on earlier drafts of this paper. Finally, I thank the anonymous referee for a careful reading.

\section{Commutative algebra}

\subsection{The height-amplitude theorem}

Let $R$ be a local Noetherian ring with maximal ideal $\mathfrak{m}$
and residue field $k$, and let $\mathbf{D}_{\mathrm{fg}}^{-}(R)$
denote the derived category of bounded-above $R$-module complexes
with finitely generated cohomology in each degree. Given $C^{\bullet}\in\mathbf{D}_{\mathrm{fg}}^{-}(R)$,
we set\[
\tau^{i}(C^{\bullet})=\mathrm{dim}_{k}H^{i}(C^{\bullet}\otimes_{R}^{\mathbf{L}}k);\]
the hypertor spectral sequence shows that $H^{i}(C^{\bullet}\otimes_{R}^{\mathbf{L}}k)$,
as a $k$-vector space, is isomorphic to a direct sum of subquotients
of $\mathrm{Tor}_{j}^{R}(H^{i+j}(C^{\bullet}),k)$, from which the
finiteness of the $\tau^{i}$'s follows easily. Any complex $C^{\bullet}\in\mathbf{D}_{\mathrm{fg}}^{-}(R)$
has a unique \emph{minimal resolution: }a bounded-above complex $F^{\bullet}$
of free $R$-modules quasi-isomorphic to $C^{\bullet}$ and such that
$\mathrm{im}d_{F^{\bullet}}^{i}\subseteq\mathfrak{m}F^{i+1}$ for
all $i$. For the existence of minimal resolutions, see  \cite{Roberts}.
A simple calculation shows that $\mathrm{rank}_{R}F^{i}=\tau^{i}(C^{\bullet})$.
By Nakayama's lemma, the greatest integer $i$ such that $\tau^{i}(C^{\bullet})\neq0$
coincides with the greatest integer $j$ such that $H^{j}(C^{\bullet})\neq0$;
we denote their common value by $d^{+}(C^{\bullet})$ or simply by
$d^{+}$ if $C^{\bullet}$ is clear.

Given $C^{\bullet}\in\mathbf{D}_{\mathrm{fg}}^{-}(R)$, we define
the \emph{amplitude }of $C^{\bullet}$ as the difference\[
\mathrm{am}(C^{\bullet})=\mathrm{sup}\left\{ i|\tau^{i}(C^{\bullet})\neq0\right\} -\mathrm{inf}\left\{ i|\tau^{i}(C^{\bullet})\neq0\right\} .\]
In general the amplitude need not be finite: if $M$ is a finite $R$-module,
viewed as a complex concentrated in degree zero, then $\mathrm{am}(M)=\mathrm{projdim}(M)$.
Note that the amplitude is finite if and only if the minimal resolution
of $C^{\bullet}$ is a bounded complex.

The first two parts of the following theorem and their proofs are
implicit in James Newton's appendix to  \cite{Universal}.

\textbf{Theorem 2.1.1. }\emph{Suppose $R$ is Cohen-Macaulay and $C^{\bullet}\in\mathbf{D}_{\mathrm{fg}}^{-}(R)$
is a complex of finite amplitude.}
\begin{description}
\item [{i.}] \emph{Any minimal prime $\mathfrak{p}$ in the $R$-support
of $H^{\ast}(C^{\bullet})$ satisfies\[
\mathrm{ht}\, \mathfrak{p}\leq\mathrm{am}(C^{\bullet}).\]
}
\item [{ii.}] \emph{If $\mathfrak{p}$ is a minimal prime in the $R$-support
of $H^{\ast}(C^{\bullet})$ with $\mathrm{ht}\, \mathfrak{p}=\mathrm{am}(C^{\bullet})$,
then $H^{j}(C^{\bullet})_{\mathfrak{p}}=0$ for $j\neq d^{+}$.}
\item [{iii.}] \emph{If $\mathrm{ht}\,  \mathfrak{p}=\mathrm{am}(C^{\bullet})$
for }every \emph{minimal prime in the $R$-support of $H^{\ast}(C^{\bullet})$,
then $H^{j}(C^{\bullet})=0$ for $j\neq d^{+}$, and $H^{d^{+}}(C^{\bullet})$
is a perfect $R$-module.}
\end{description}
\emph{Proof of i. and ii. }Replacing $C^{\bullet}$ by its minimal
resolution, we may assume $C^{\bullet}$ is a bounded complex of free
$R$-modules of finite rank (and as such, we may write derived tensor
products of $C^{\bullet}$ as ordinary tensor products). Let $d^{-}$
be the least integer $i$ for which $\tau^{i}(C^{\bullet})\neq0$.
$ $Let $\mathfrak{p}$ be a minimal element of $\mathrm{Supp}H^{\ast}(C^{\bullet})$,
and let $r$ be the least degree with $\mathfrak{p}\in\mathrm{Supp}H^{r}(C^{\bullet})$.
Let $h=\mathrm{ht}\, \mathfrak{p}$, and choose a system of parameters
$x_{1},\dots,x_{h}\in\mathfrak{p}$ for $R_{\mathfrak{p}}$. Set $J_{n}=(x_{1},\dots,x_{n})$.
We will show inductively that $H^{r-n}(C^{\bullet}\otimes R/J_{n})\neq0$
for $1\leq n\leq h$. Granted this inductive step, the theorem follows
from the following observation: letting $\check{C}_{\bullet}=\mathrm{Hom}_{R}(C^{\bullet},R)$
denote the dual complex, there is a natural spectral sequence\[
E_{2}^{i,j}=\mathrm{Ext}_{R}^{i}(H_{j}(\check{C}_{\bullet}),R/I)\Rightarrow H^{i+j}(C^{\bullet}\otimes R/I).\]
Since $H^{r-h}(C^{\bullet}\otimes R/J_{h})\neq0$, the least $j$
with $H_{j}(\check{C}_{\bullet})\neq0$, say $j_{\mathrm{min}}$,
satisfies $j_{\mathrm{min}}\leq r-h$. Taking $I=\mathfrak{m}$, the
entry $E_{2}^{0,j_{\mathrm{min}}}$ is stable and nonzero, so $d^{-}=j_{\mathrm{min}}\leq r-h$.
Putting things together gives\[
d^{-}+h\leq r\leq d^{+},\]
so $h\leq d_{\mathfrak{}}^{+}-d^{-}=\mathrm{am}(C^{\bullet})$, as
desired. If equality holds then $r\geq d^{-}+\mathrm{am}(C^{\bullet})=d^{+}$,
but $r\leq d^{+}$ was the \emph{least }degree with $\mathfrak{p}\in\mathrm{Supp}H^{r}(C^{\bullet})$.

It remains to carry out the inductive step. Let $\mathfrak{p}_{n}$
denote the image of $\mathfrak{p}$ under $R_{\mathfrak{p}}\to R_{\mathfrak{p}}/J_{n}$.
For $0\leq n\leq h-1$, suppose $H^{r-n}(C^{\bullet}\otimes R/J_{n})_{\mathfrak{p}_{n}}$
is nonzero with $\mathfrak{p}_{n}$ an associated prime, and $H^{i}(C^{\bullet}\otimes R/J_{n})_{\mathfrak{p}_{n}}=0$
for $i<r-n$. Then $H^{r-n-1}(C^{\bullet}\otimes R/J_{n+1})_{\mathfrak{p}_{n+1}}$
is nonzero with $\mathfrak{p}_{n+1}$ an associated prime, and $H^{i}(C^{\bullet}\otimes R)_{\mathfrak{p}_{n+1}}=0$
for $i<r-n-1$. The supposition is true for $n=0$ by our assumptions
and the fact that minimal primes are associated primes. To prove the
induction, we proceed as follows. For each $0\leq n\le h-1$ we have
a spectral sequence\[
E_{2}^{i,j}=\mathrm{Tor}_{-i}^{R/J_{n}}(H^{j}(C^{\bullet}\otimes R/J_{n}),R/J_{n+1})\Rightarrow H^{i+j}(C^{\bullet}\otimes R/J_{n+1})\]
of $R/J_{n}$-modules. Localize this spectral sequence at $\mathfrak{p}$.
Since $R_{\mathfrak{p}}$ is Cohen-Macaulay, any system of parameters
is a regular sequence on $R_{\mathfrak{p}}$. As such, caculating
$\mathrm{Tor}^{R_{\mathfrak{p}}/J_{n}}(-,R_{\mathfrak{p}}/J_{n+1})$
via the resolution \[
0\to R_{\mathfrak{p}}/J_{n}\overset{\cdot x_{n+1}}{\to}R_{\mathfrak{p}}/J_{n}\to R_{\mathfrak{p}}/J_{n+1}\to0\]
implies that the entries of the spectral sequence vanish for $i\neq0,1$,
with $E_{2}^{-1,j}\simeq H^{j}(C^{\bullet}\otimes R/J_{n})_{\mathfrak{p}_{n}}[x_{n+1}]$.
The vanishing claim follows easily, and we get an isomorphism \[
H^{r-n-1}(C^{\bullet}\otimes R/J_{n+1})_{\mathfrak{p}_{n+1}}\simeq H^{r-n}(C^{\bullet}\otimes R/J_{n})_{\mathfrak{p}_{n}}[x_{n+1}],\]
of $R_{\mathfrak{p}}/J_{n+1}$-modules; by our inductive hypothesis
the right-hand side is easily seen to be nonzero with $\mathfrak{p}_{n+1}$
an associated prime. This completes the proof of i. and ii.

\emph{Proof of iii.} Let\[
F^{\bullet}:0\to F^{0}\to F^{1}\to\dots\to F^{d}\to0\]
be a complex of free finite rank $R$-modules such that every minimal
prime in the $R$-support of $H^{\ast}(F^{\bullet})$ has height exactly
$d$. By parts i. and ii., every minimal prime in the support of $H^{i}(F^{\bullet})$
for $0\leq i\leq d-1$ has height $\geq d+1$. Consider the dual complex
$\check{F}^{\bullet}=\mathrm{Hom}_{R}(F^{-\bullet},R)$. A priori
the cohomology of $\check{F}^{\bullet}$ is concentrated in degrees
$-d$ through $0$, and we have a convergent spectral sequence \[
E_{2}^{i,j}=\mathrm{Ext}_{R}^{i}(H^{-j}(F^{\bullet}),R)\Rightarrow H^{i+j}(\check{F}^{\bullet}).\]
Since $R$ is Cohen-Macaulay, we have $\mathrm{grade}M+\mathrm{dim}M=\mathrm{dim}R$
for any $R$-module $M$, and thus the entries $E_{2}^{i,j}$ vanish
for $j<-d$, for $j=-d$ with $i<d$, and for $j>-d$ with $i<d+1$.
Thus the spectral sequence yields isomorphisms $H^{i}(\check{F}^{\bullet})=0$
for $i\neq0$ and $H^{0}(\check{F}^{\bullet})\simeq\mathrm{Ext}_{R}^{d}(H^{d}(F^{\bullet}),R)$.
Applying the adjunction isomorphism\[
\mathbf{R}\mathrm{Hom}_{R}(\mathbf{R}\mathrm{Hom}_{R}(F^{\bullet},R),R)\cong F^{\bullet}\]
yields a dual spectral sequence\[
E_{2}^{i,j}=\mathrm{Ext}_{R}^{i}(H^{-j}(\check{F}^{\bullet}),R)\Rightarrow H^{i+j}(F^{\bullet}),\]
which correspondingly degenerates to an isomorphism $H^{i}(F^{\bullet})\simeq\mathrm{Ext}_{R}^{i}(H^{0}(\check{F}^{\bullet}),R)$
for any $i$. Since $\mathrm{Ann}H^{0}(\check{F}^{\bullet})\supseteq\mathrm{Ann}H^{d}(F^{\bullet})$,
we have $\mathrm{grade}H^{0}(\check{F}^{\bullet})\geq\mathrm{grade}H^{d}(F^{\bullet})=d$,
so $H^{i}(F^{\bullet})$ vanishes for $i<d$. Thus $F^{\bullet}$
yields a free resolution of $H^{d}(F^{\bullet})$, so $\mathrm{projdim}H^{d}(F^{\bullet})\leq d$.
Quite generally we have $\mathrm{grade}M\leq\mathrm{projdim}M$, so
perfection follows. $\square$

\subsection{A patching theorem for complexes}

Fix a complete discrete valuation ring $\mathcal{O}$. Set $R_{\infty}=\mathcal{O}[[x_{1},\dots,x_{q-r}]]$
and $S_{\infty}=\mathcal{O}[[T_{1},\dots,T_{q}]]$. Write $S_{n}$
for the quotient $S_{\infty}/\left((1+T_{1})^{p^{n}}-1,\dots,(1+T_{q})^{p^{n}}-1\right)$,
with $\overline{S}_{n}=S_{n}/\varpi^{n}$. We write $\mathfrak{a}$
for the ideal $(T_{1},\dots,T_{q})$ in $S_{\infty}$ and in $S_{n}$,
and we abusively write $k$ for the common residue field of all these
local rings.

\textbf{Theorem 2.2.1. }\emph{Let $R$ be a complete local Noetherian
$\mathcal{O}$-algebra, and let $H$ be an $R$-module which is $\mathcal{O}$-module
finite. Suppose for each integer $n\geq1$ we have a surjection $\phi_{n}:R_{\infty}\to R$
and a complex $C_{n}^{\bullet}\in\mathbf{D}_{\mathrm{fg}}^{-}(S_{n})$
of $S_{n}$-modules with the following properties:}
\begin{description}
\item [{i.}] \emph{$\tau^{i}(C_{n}^{\bullet})$ is independent of $n$
for $i\in[d,d-r]$ and zero for $i\notin[d,d-r]$, where $d$ is some
fixed integer.}
\item [{ii.}] \emph{There is a degree-preseving $R_{\infty}$-module structure
on $H_{n}^{\ast}=H^{\ast}(C_{n}^{\bullet})$ such that the image of
$S_{n}$ in $\mathrm{End}_{\mathcal{O}}(H_{n}^{\ast})$ arises from
an algebra homomorphism $i_{n}:S_{\infty}\to R_{\infty}$ with $(\phi_{n}\circ i_{n})(\mathfrak{a})=0$.}
\item [{iii.}] \emph{Writing $H_{n}=H^{d}(C_{n}^{\bullet})$, there is
an isomorphism $H_{n}/\mathfrak{a}H_{n}\simeq H$ of $R_{\infty}$-modules
where $ $$R_{\infty}$ acts on $H$ via $\phi_{n}$.}
\end{description}
\emph{Then $H$ is free over $R$.}

\emph{Proof. }Let $F_{n}^{\bullet}$ be the minimal resolution of
$C_{n}^{\bullet}$, and set $\overline{F}_{n,m}^{\bullet}=F_{n}^{\bullet}\otimes_{S_{n}}\overline{S}_{m}$
for $m\leq n$. By our assumptions $\mathrm{dim}_{\overline{S}_{m}}\overline{F}_{n,m}^{i}=\tau^{i}(C_{n}^{\bullet})=\tau^{i}$
is independent of $n$ and $m$. Choosing bases we can and do represent
the differentials $d_{n,m}^{i}$ of $\overline{F}_{n,m}^{\bullet}$
by matrices $\delta_{n,m}^{i}\in M_{\tau^{i}\times\tau^{i+1}}(\overline{S}_{m})$.
By the usual pigeonhole argument we may find integers $j_{n},n\geq1$
such that $\delta_{j_{n},n}^{i}$ is the reduction of $\delta_{j_{n+1},n+1}^{i}$.
Let $\delta_{\infty}^{i}$ be the limit of the sequence $\delta_{j_{n},n}^{i}$
as $n\to\infty$, and let $F_{\infty}^{\bullet}$ be the bounded complex
of free finite rank $S_{\infty}$-modules whose differentials are
given by the $\delta_{\infty}^{i}$'s. Set $H_{\infty}^{\ast}=H^{\ast}(F_{\infty}^{\bullet})$
and $H_{\infty}=H^{d}(F_{\infty}^{\bullet})$. Passing to a further
subsequence if necessary, the maps $i_{n}$ and $\phi_{n}$ converge
to algebra homomorphisms $i_{\infty}$ and $\phi_{\infty}$, and the
$R_{\infty}$-module structures on $H_{j_{n}}^{\ast}$ patch together
into an $R_{\infty}$-module structure on $H_{\infty}^{\ast}$ such
that $S_{\infty}$ acts through $i_{\infty}$. Since $H_{\infty}^{\ast}$
is finite over $S_{\infty}$, it is finite over $R_{\infty}$, so
in particular \[
\mathrm{dim}_{S_{\infty}}H_{\infty}^{\ast}=\mathrm{dim}_{R_{\infty}}H_{\infty}^{\ast}\leq\mathrm{dim}R_{\infty}=1+q-r.\]
On the other hand, the first part of the height-amplitude theorem
implies the \emph{opposite }inequality, so every minimal prime in
the $S_{\infty}$-support of $H_{\infty}^{\ast}$ has height exactly
$r$. Therefore, $H_{\infty}\simeq H_{\infty}^{\ast}$ by the third
part of the height-amplitude theorem, and $F_{\infty}^{\bullet}$
is a free resolution of $H_{\infty}$ of length $r$. This shows that
$\mathrm{projdim}_{S_{\infty}}(H_{\infty})=r$, so $\mathrm{depth}_{S_{\infty}}(H_{\infty})=1+q-r$
by the Auslander-Buchsbaum formula. But then $\mathrm{depth}_{R_{\infty}}(H_{\infty})=1+q-r=\mathrm{dim}R_{\infty}$
via $i_{\infty}$, so $H_{\infty}$ is a free module over $R_{\infty}$
by a second application of Auslander-Buchsbaum. Therefore $H_{\infty}/\mathfrak{a}H_{\infty}$
is a free module over $R_{\infty}/i_{\infty}(\mathfrak{a})$. But
$H_{\infty}/\mathfrak{a}H_{\infty}\simeq H$ as $R_{\infty}/i_{\infty}(\mathfrak{a})$-modules,
where $R_{\infty}/i_{\infty}(\mathfrak{a})$ acts on $H$ through
the surjection $R_{\infty}/i_{\infty}(\mathfrak{a})\twoheadrightarrow R$
induced by $\phi_{\infty}$. $\square$

\section{Modularity lifting}

We return to the notation of the introduction. Let $S(\overline{\rho})$
be the ramification set of $\overline{\rho}$, let $Q$ be any finite
set of primes disjoint from $S(\overline{\rho})\cup\{v|p\}$, and
let $S_{Q}$ denote the set of places $Q\cup S(\overline{\rho})\cup\{v|p\}$.
For any such $Q$, let $R_{Q}$ denote the deformation ring defined
in §4.1 of  \cite{CG}; this is a complete local Noetherian $\mathcal{O}$-algebra.
Let $H_{Q}^{1}(F,\mathrm{ad}^{0}\overline{\rho})$ be the Selmer group
defined as the kernel of the map\[
H^{1}(F,\mathrm{ad}^{0}\overline{\rho})\to\prod_{v}H^{1}(F_{v},\mathrm{ad}^{0}\overline{\rho})/L_{v}\]
where $L_{v}=H_{\mathrm{ur}}^{1}(F_{v},\mathrm{ad}^{0}\overline{\rho})$
if $v\notin Q\cup\{v|p\}$, $L_{v}=H^{1}(F_{v},\mathrm{ad}^{0}\overline{\rho})$
if $v\in Q$, and $L_{v}=H_{\mathrm{f}}^{1}(F_{v},\mathrm{ad}^{0}\overline{\rho})$
if $v|p$ (here $H_{\mathrm{f}}^{1}$ is as in §2.4 of  \cite{DDT}).
Modifying the proof of Corollary 2.43 of  \cite{DDT} via Corollary 2.4.3 of \cite{CHT}, we find that
the reduced tangent space of $R_{Q}$ has dimension at most \[
\mathrm{dim}_{k}H_{Q}^{1}(F,\mathrm{ad}^{0}\overline{\rho}(1))-r+\sum_{v\in Q}\mathrm{dim}_{k}H^{0}(F_{v},\mathrm{ad}^{0}\overline{\rho}(1)).\]

We define $L_{Q}$ and $K_{Q}$ be the open compact subgroups as in
 \cite{CG}. We denote by $Y_{0}(Q)$ the arithmetic quotient $\mathrm{GL}_{2}(F)\backslash\mathrm{GL}_{2}(\mathbf{A}_{F})/L_{Q}K_{\infty}^{\circ}Z_{\infty}$,
and by $Y_{1}(Q)$ the quotient $\mathrm{GL}_{2}(F)\backslash\mathrm{GL}_{2}(\mathbf{A}_{F})/K_{Q}K_{\infty}^{\circ}Z_{\infty}$.
For $Q=\emptyset$ we simply write $Y$. Let \[
T_{v}=L_{Q}\left(\begin{array}{cc}
\varpi_{v}\\
 & 1\end{array}\right)L_{Q}\]
and\[
\left\langle v\right\rangle =L_{Q}\left(\begin{array}{cc}
\varpi_{v}\\
 & \varpi_{v}\end{array}\right)L_{Q}\]
denote the usual Hecke operators; when $v\in Q$ we write $U_{v}$
for $T_{v}$ as is customary. Let $\mathscr{T}_{Q}^{\mathrm{an}}$
denote the abstract polynomial algebra over $\mathcal{O}$ in the
operators $T_{v}$ and $\left\langle v\right\rangle $ for all
places $v\notin S_{Q}$, and let $\mathscr{T}_{Q}$ denote the algebra
generated by these operators together with the operators $U_{v}$
for $v\in Q$. Write $\mathbf{T}_{Q}^{\mathrm{an}}$ and $\mathbf{T}_{Q}$
for the images of $\mathscr{T}_{Q}^{\mathrm{an}}$ and $\mathscr{T}_{Q}$
in $\mathrm{End}_{\mathcal{O}}\left(H^{\ast}(Y_{1}(Q),\mathcal{O})\right)$.
When $Q=\emptyset$ we write $\mathbf{T}=\mathbf{T}_{\emptyset}$.
By assumption $\overline{\rho}$ is associated with a maximal ideal
of $\mathbf{T}$ which we denote by $\mathfrak{m}_{\emptyset}$.

Let $\mathfrak{m}$ be any maximal ideal of $\mathbf{T}_{Q}$ containing
the preimage of $\mathfrak{m}_{\emptyset}$ under the natural map
$\mathbf{T}_{Q}^{\mathrm{an}}\to\mathbf{T}$. We make the following
conjecture:

\textbf{Conjecture 3.1 (Existence Conjecture): }\emph{For any $Q$
there is a surjection $\phi_{Q}:R_{Q}\twoheadrightarrow\mathbf{T}_{Q,\mathfrak{m}}$
such that the associated Galois representation $\rho_{Q}:G_{F}\to\mathrm{GL}_{2}(\mathbf{T}_{Q,\mathfrak{m}})$
has the following properties:}
\begin{description}
\item [{i.}] \emph{For any $v\notin S_{Q}$,\[
\det\left(X-\rho_{Q}(\mathrm{Frob}_{v})\right)=X^{2}-T_{v}X+\mathbf{N}v\cdot\left\langle v\right\rangle \in\mathbf{T}_{Q,\mathfrak{m}}[X].\]
}
\item [{ii.}] \emph{For any $v\in Q$, $\rho_{Q}|D_{v}\simeq\eta_{1}\oplus\eta_{2}$
with $\eta_{1}$ unramified and $\eta_{1}(\mathrm{Frob}_{v})=U_{v}$.}
\end{description}
\emph{Furthermore, $H^{i}(Y_{1}(Q),\mathcal{O})\otimes_{\mathbf{T}_{Q}}\mathbf{T}_{Q,\mathfrak{m}}$
vanishes for $i\notin[d-r,d]$, where $d=[F:\mathbf{Q}]$.}

Suppose now for each $n\geq1$ that $Q_{n}$ is a set of Taylor-Wiles
primes of cardinality $q=\mathrm{dim}_{k}H_{\emptyset}^{1}(F,\mathrm{ad}^{0}\overline{\rho}(1))$,
such that each $v\in Q_{n}$ has $\mathbf{N}v\equiv1\,\mathrm{mod}\, p^{n}$.
The reduced tangent space of $R_{Q_{n}}$ has dimension at most $q-r$.
Let $\mathfrak{m}_{n}$ denote the maximal ideal of $\mathbf{T}_{Q_{n}}$
generated by the preimage of $\mathfrak{m}$ under the map $\mathbf{T}_{Q_{n}}^{\mathrm{an}}\to\mathbf{T}$
and by $U_{v}-\alpha_{v}$ for all $v\in Q_{n}$, where $\alpha_{v}$
is a fixed choice of one of the eigenvalues of $\overline{\rho}(\mathrm{Frob}_{v})$.

\textbf{Proposition 3.2. }\emph{There is an isomorphism\[
H^{\ast}(Y_{0}(Q_{n}),\mathcal{O})_{\mathfrak{m}_{n}}\simeq H^{\ast}(Y,\mathcal{O})_{\mathfrak{m}}.\]
}This is proved exactly as in Lemmas 3.4 and 4.6 of  \cite{CG}, working
one degree at a time.

By design there is a natural surjection \[
\prod_{v\in Q_{n}}(\mathcal{O}_{F_{v}}/\varpi_{v})^{\times}\twoheadrightarrow(\mathbf{Z}/p^{n})^{q}.\]
Composing this with the natural reduction map $L_{Q_{n}}\to\prod_{v\in Q_{n}}(\mathcal{O}_{F_{v}}/\varpi_{v})^{\times}$
gives $S_{n}=\mathcal{O}[(\mathbf{Z}/p^{n})^{q}]$ the structure of
a local system over $Y_{0}(Q_{n})$. Let $X_{\mathbf{A}}$ be the
quotient $\mathrm{GL}_{2}(F)\backslash\mathrm{GL}_{2}(\mathbf{A}_{F})/K_{\infty}^{\circ}Z_{\infty}$,
and let $C_{\bullet}(X_{\mathbf{A}})$ be the complex of singular
chains on $X_{\mathbf{A}}$ with $\mathbf{Z}$-coefficients. Set $C^{\bullet}(Q_{n})=\mathrm{Hom}_{\mathbf{Z}}(C_{\bullet}(X_{\mathbf{A}}),S_{n})^{L_{Q_{n}}}$,
so there is a canonical isomorphism $H^{\ast}(C^{\bullet}(Q_{n}))\simeq H^{\ast}(Y_{0}(Q_{n}),S_{n})$.
The canonical action of $\mathrm{GL}_{2}(\mathbf{A}_{F}^{f})$ on
$C_{\bullet}(X_{\mathbf{A}})$ induces a canonical action of the algebra $\mathscr{T}_{Q_{n}}$
on the complex $C^{\bullet}(Q_{n})$ lifting the Hecke action on cohomology: precisely, given $\phi \in \mathrm{Hom}_{\mathbf{Z}}(C_{\bullet}(X_{\mathbf{A}}),S_{n})^{L_{Q_{n}}}$ regarded as a function on chains, $T_g=[L_{Q_n}gL_{Q_n}] \in \mathscr{T}_{Q_n}$ acts on $\phi$ by \[(T_g \cdot \phi)(\sigma)=\sum_{g_i \in L_{Q_n}gL_{Q_n}/L_{Q_n}} g_i\cdot \phi(\sigma g_i).\] (For a more thorough discussion of this idea, which the author learned from Glenn Stevens, see \cite{Universal}.)

Let $\mathfrak{M}_{n}$ be the
preimage of $\mathfrak{m}_{n}$ under the structure map $\mathscr{T}_{Q_{n}}\twoheadrightarrow\mathbf{T}_{Q_{n}}$,
and set \[
C_{n}^{\bullet}=C^{\bullet}(Q_{n})\otimes_{\mathscr{T}_{Q_{n}}}\mathscr{T}_{Q_{n},\mathfrak{M}_{n}}.\]
Since $\mathscr{T}_{Q_{n},\mathfrak{M}_{n}}$ is flat over $\mathscr{T}_{Q_{n}}$, the functor
$-\otimes_{\mathscr{T}_{Q_{n}}}\mathscr{T}_{Q_{n},\mathfrak{M}_{n}}$
commutes with taking cohomology, so we have canonical isomorphisms
\begin{eqnarray*}
H^{\ast}(C_{n}^{\bullet}) & \simeq & H^{\ast}(C^{\bullet}(Q_{n}))\otimes_{\mathscr{T}_{Q_{n}}}\mathscr{T}_{Q_{n},\mathfrak{M}_{n}}\\
 & \simeq & H^{\ast}(Y_{0}(Q_{n}),S_{n})\otimes_{\mathbf{T}_{Q_{n}}}\mathbf{T}_{Q_{n}}\otimes_{\mathscr{T}_{Q_{n}}}\mathscr{T}_{Q_{n},\mathfrak{M}_{n}}\\
 & \simeq & H^{\ast}(Y_{0}(Q_{n}),S_{n})_{\mathfrak{m}_{n}}.\end{eqnarray*}
If $\mathfrak{a}$ denotes the augmentation ideal of $S_{n}$
then\begin{eqnarray*}
H^{i}(C_{n}^{\bullet}\otimes_{S_{n}}^{\mathbf{L}}S_{n}/\mathfrak{a}S_{n}) & \simeq & H^{i}\left(Y_{0}(Q_{n}),\mathcal{O}\right)_{\mathfrak{m}_{n}}\\
 & \simeq & H^{i}(Y,\mathcal{O})_{\mathfrak{m}}\end{eqnarray*}
by Proposition 3.2. This shows that the complexes $C_{n}^{\bullet}$
satisfy assumption i. of Theorem 2.2.1. We are ready to verify the
rest of the assumptions of Theorem 2.2.1. For each $n$, fix a choice
of a surjection $\sigma_{n}:R_{\infty}\twoheadrightarrow R_{Q_{n}}$.
The composite $\phi_{Q_{n}}\circ\sigma_{n}$ where $\phi_{Q_{n}}$
is the map provided by the existence conjecture gives $H^{\ast}(C_{n}^{\bullet})$
a degree-preserving $R_{\infty}$-module structure. Define $\phi_{n}$
as the composite of our chosen surjection $R_{\infty}\twoheadrightarrow R_{Q_{n}}$
with the natural surjection $R_{Q_{n}}\twoheadrightarrow R^{\mathrm{min}}$.
As in §2.8 of  \cite{DDT}, $R_{Q_{n}}$ is naturally an $S_{n}$-algebra,
with $R_{Q_{n}}/\mathfrak{a}R_{Q_{n}}\simeq R^{\mathrm{min}}$, and
the map $\phi_{Q_{n}}$ is equivariant for the $S_{n}$-actions on
its source and target. Let $i_{n}:S_{\infty}\to R_{\infty}$ be any
fixed lift of the composite $S_{\infty}\to S_{n}\to R_{Q_{n}}$ compatible
with $\sigma_{n}$, so $(\phi_{n}\circ i_{n})(\mathfrak{a})=0$ by
construction. We've now verified assumption ii. Assumption iii. is
immediate from Hochschild-Serre, so Theorem 2.2.1 applies and modularity
lifting at minimal level follows.

\bibliographystyle{alpha}

\end{document}